\documentclass[11pt]{extarticle}

\usepackage[left=1in,right=1in,top=1in,bottom=1.3in]{geometry}
\usepackage{amsmath, amssymb, amsthm, latexsym, amsfonts, indentfirst, xcolor, mathtools, manfnt, microtype}
\usepackage[utf8]{inputenc}
\usepackage[english]{babel}

\numberwithin{equation}{section} % нумерация формул по разделам

\usepackage{indentfirst}
\usepackage{booktabs} 		% таблицы

\usepackage{tikz}
\usetikzlibrary{decorations.pathreplacing,calligraphy}
\usetikzlibrary{shapes.geometric}

\usepackage{amsthm}
\theoremstyle{plain}
\newtheorem{theorem}{Theorem}
\newtheorem{proposition}{Proposition}
\newtheorem{lemma}{Lemma}
\newtheorem{corollary}{Corollary}

\theoremstyle{definition}

\DeclareMathOperator{\E}{\mathrm E}
\newcommand{\ff}{\mathcal{F}}
\newcommand{\G}{\mathcal{G}}
\newcommand{\FF}{\mathcal{F}}
\newcommand{\m}{\mathcal}
\usepackage{hyperref}

\newenvironment{proofw}{\par
  \pushQED{\qed}%
  \normalfont
  \trivlist
  \item[]\ignorespaces
}{%
  \popQED\endtrivlist\
}

\title{Two questions on Kneser colorings}
\date{}
\author{Eduard Inozemtsev\footnote{Moscow Institute of Physics and Technology, Email: {\tt eduard\_inozemtsev@bk.ru}}, Andrey Kupavskii\footnote{Moscow Institute of Physics and Technology, Saint-Petersburg State University Russia; Email: {\tt kupavskii@ya.ru}.
The research was supported by the grant of the Russian Science Foundation N 24-71-10021.
}
}
\usepackage{enumitem}
\begin{document}
\maketitle
\begin{abstract}
    In this paper, we investigate two questions on Kneser graphs $KG_{n,k}$. First, we prove that the union of $s$ intersecting families in ${[n]\choose k}$ has size at most ${n\choose k}-{n-s\choose k}$ for all sufficiently large $n$ that satisfy $n>(2+\epsilon)k^2+s$ with $\epsilon>0$. We provide an example that shows that this result is essentially tight for the number of colors close to $\chi(KG_{n,k})=n-2k+2$. We also improve the result of Bulankina and Kupavskii on the choice chromatic number, showing that it is at least $\frac 1{25} n\log n$ for all $k<\sqrt n$ and $n$ sufficiently large.
\end{abstract}
\section{Introduction}

For positive integers $n, k$, let $[n]:=\{ 1, \dots, n\}$ denote the standard $n$-element set and let ${[n] \choose k} := \{ F \subset [n] : |F| = k \}$ be the set of its $k$-element subsets.  For a set $X$ we denote its power set by $2^X$.  A {\it family} is a collection of sets  $\FF \subset 2^X$. We say that a family $\FF$ is \emph{intersecting} if for any $A, B \subset \mathcal{F}$, we have $A \cap B \neq \emptyset$. One of the classical results in extremal combinatorics is the Erd\H{o}s--Ko--Rado theorem \cite{EKR}, which states that for $n\ge 2k$ the largest intersecting family in ${[n]\choose k}$ has size at most ${n-1\choose k-1}$.  For $n > 2k$ equality holds only if $\FF$ is of the form $\mathcal{S}_x := \{ F \subset {[n] \choose k} : x \in F \}$, for some $x \in [n]$. The {\it covering number} $\tau(\ff)$ of $\ff$ is the size of the smallest set $X$ such that $X\cap F\ne \emptyset$ for all $F\in \ff$. Then $\mathcal{S}_x$ has covering number $1$. Such families are called \emph{stars} or \emph{trivial}, and $x$ is called their {\it center}. Hilton and Milner \cite{HM} gave a tight upper bound on the size of intersecting family in ${[n]\choose k}$ with covering number at least $2$ for $n>2k$:
\begin{equation}\label{eqhm}
    |\ff|\le {n-1\choose k-1}-{n-k-1\choose k-1}+1.
\end{equation}
The extremal families have the property that all but one `exceptional' set can be pierced by a single element.
%the following tight result about the size of non-trivial intersecting families:
%\begin{theorem}[Hilton--Milner]\label{HM}
%    Let $n > 2k$. If $\FF$ is a non-trivially intersecting family, then $|\FF| \leq {n-1 \choose k-1} - {n-k-1 \choose k-1} + 1$. Equality holds only if $\FF = \{ F \subset {[n] \choose k} : x \in F, F \cap  H \neq \emptyset \} \cap H$ , for some $x \in [n]$ and $H \subset {[n] \choose k}.$
%\end{theorem}

The following generalization of the Erd\H{o}s--Ko--Rado theorem was  proposed by Erd\H{o}s \cite{ERD}: for given $n,k,s$, determine, what is the largest size of a family $\FF \subset {[n] \choose k}$ such that $\FF$ is a union of at most $s$ intersecting families. It is a natural conjecture to make  that for triples of parameters $(n, k, s)$ where $n$ is somewhat large in comparison to $k,s$, the extremal example is a union of $s$ stars. Erd\H os indeed conjectured this for all $n\ge 2k+s-1$. In this case, $|\FF| = {n \choose k} - {n-s \choose k}$. Frankl and F\" uredi \cite{FF} showed that the conjecture is true for $s = 2$ and $ n > ck$ with $c \approx 2.62$. However, it is not true in general. This was shown  for $n=2k+s-1$ by Hilton and F\"uredi. Frankl and F\"uredi also gave an example which shows that the extreme is different even for $s=2$, when $n \leq 2k + c_0 \sqrt{k}$ for $n\ge n_0$ and some positive constant $c_0$.  Ellis and Lifshitz \cite{EL} recently proved that for constant $s$ and $n > 2k + C(s) k^{2/3}$ the union of stars is optimal, where  $C(s)$ is a function of $s$.

It is very natural to interpret this problem in terms of Kneser graphs. Consider a \emph{Kneser graph} $KG_{n,k}$,  whose vertex set is ${[n]\choose k}$, and edges connect pairs of disjoint sets. Kneser showed that $\chi(KG_{n,k})\le n-2k+2$ and conjectured that this is tight. Lov\'asz proved \cite{L}  that  this indeed holds. By now, there is a series of results on the chromatic number of Kneser graphs and different subgraphs (e.g., see the recent paper of Kaiser and Stehl\'ik \cite{KS}, where they construct edge-critical subgraphs of Kneser graphs with the same chromatic number). However, all these proofs are topological and thus are difficult to adapt to other similar settings. The relation to the above question is not difficult to discern. Actually, intersecting families in ${[n]\choose k}$ are exactly the independent sets in $KG_{n,k},$ and thus are the collections of vertices that can be colored by the same color. In these terms, the question of Erd\H os may be reformulated as follows: what is the maximum number of vertices of $KG_{n,k}$ that can be properly colored in $s$ colors?

In this paper, we answer this question in the regime when the number of colors can be very large, but $n$ is somewhat large compared to $k$.

\begin{theorem}\label{main}
Let $\epsilon>0$ and let $n,k$ be integers such that $n>(2+\epsilon)k^2$ and $n$ is sufficiently large. If $\FF$ is a union of $s$ non-trivial intersecting families, $s<n$, then $|\FF| <  {n \choose k} - {n-s \choose k}$.
\end{theorem}

A simple argument using the Hilton--Milner theorem gives that the conclusion of the theorem is valid for, roughly speaking, $n\ge 2k^3$. We present it in the next section.

Why do we consider the case when all families are non-trivial?  If $\FF_1, \dots, \FF_s \subset 2^{[n]}$ are $s$ different stars, then clearly $|\cup_{i \in [s]} \FF_i| = {n \choose k} - {n-s \choose k}$. Now let $\FF_1$ be a star with the center in $1$. % and let $\FF_2, \dots, \FF_s$ be non-trivially intersecting.
Remove $1$ from the ground set, also removing all sets which contain $1$, and $\ff_1$ in particular. The problem for 
$(n, k, s)$ is reduced to the problem for  $(n-1, k, s-1)$. Assuming that we have $|\cup_{i=2}^s \ff_i|\le {n-1\choose k}-{n-1-(s-1)\choose k} = {n-1\choose k}-{n-s\choose k}$, we get $|\cup_{i\in[s]}\ff_i|\le {n-1\choose k-1}+{n-1\choose k}-{n-s\choose k}={n\choose k}-{n-s\choose k}$. If we get rid of the assumption of all families being non-trivial then we get the following corollary.

\begin{theorem}
Let $\epsilon>0$ and let $n,s,k$ be integers such that $n>(2+\epsilon)k^2+s$ and $n$ is sufficiently large. If $\FF$ is a union of $s$ intersecting families then $|\FF| \le {n \choose k} - {n-s \choose k}$. Equality is possible only if $\ff$ is a union of $s$ stars.
  
\end{theorem}

We also show that the bounds in Theorem~\ref{main} are not far from being tight.
\begin{theorem}\label{main3}
  Let $n = (k-1)(k+1-m)+m$, where $k,m$ are positive integers and $m\le \frac k4$ and $k$ is sufficiently large. Then there exists a collection of $s:=n-2k+2-m$ non-trivial intersecting families in ${[n]\choose k}$ such that their union has size bigger than ${n\choose k}-{n-s\choose k}$.    
\end{theorem}

For the proof of Theorem~\ref{main} we rely on the techniques from the paper of Kiselev and Kupavskii \cite{KK}. There, the authors studied colorings of $KG_{n,k}$ into $n-2k+2$ colors without trivial colors, and managed to show that such colorings do not exist for $n>(2+\epsilon)k^2$.  The question of Erd\H os is related to the famous Erd\H os Matching Conjecture. In terms of $KG_{n,k}$, it asks for the largest subset of vertices with no $(s+1)$-clique. The extremal examples for $n>Csk$ are the same as for the problem discussed above, albeit the regimes are different: the EMC becomes trivial for $n<(s+1)k$. See the papers of Frankl, Kolupaev and the second author \cite{FK, KK2} for the most recent results. It is also worth mentioning a recent paper of Frankl and the second author \cite{FK2}, in which they studied the size of the largest intersecting family in ${[n]\choose k}$ with large $\tau$ and showed that there is a certain `phase transition' around $n =k^2$.

Another question that we address concerns the choice number of $KG_{n,k}.$ Recall that the choice number $ch(G)$ of a graph $G$ is the smallest $\ell$ such that, whenever each vertex of the graph is assigned a list of colors of size at least $\ell,$ we can choose one color from the list for each vertex so that the resulting coloring is proper (i.e., does not have monochromatic edges).
Bulankina and Kupavskii \cite{BK} showed that $ch(KG_{n,k})\ge \frac{1}{2\ell^2} n\log n$ for $k\le n^{\frac 12-\frac 1\ell}$ and sufficiently large $n$. They also showed the upper bound $ch(KG_{n,k})\le n\ln \frac nk +n,$ valid for any $n\ge 2k.$ In this paper, we improve their lower bound.
\begin{theorem}\label{main2}
    Fix $\varepsilon>0$ and let $n\ge n_0(\varepsilon)$. Suppose that $ 3 \leq k \leq \left(\frac{n}{1+\varepsilon}\right)^{1/2} $. Then  $$ch(KG_{n, k}) > \left(\frac 14 - \frac 1{2k}-\varepsilon \right)\frac{n\log{n}}{2}.$$
\end{theorem}
The proofs of Theorems~\ref{main} and~\ref{main2} both rely on structural properties of intersecting families that were studied in the paper of Kiselev and the second author \cite{KK}. 
\section{Proofs}
\subsection{Preliminaries}
Using the Hilton--Milner theorem, it is easy to give a bound for Theorem~\ref{main} in the following regime.
\begin{proposition} \label{hmbound}
    Let $n > k^2 + k^{3/2}\sqrt{s}$, $s\le n$, and $\FF \subset {[n] \choose k}$ is a union of $s$ non-trivial intersecting families, then $|\FF| <  {n \choose k} - {n-s \choose k}. $
\end{proposition}
We remark that it covers the case $n\ge 2k^3$. 
\begin{proof}
    The Hilton-Milner Theorem implies
    \begin{equation}\label{eqhm1}
    |\FF| \leq s\Big({n-1\choose k-1}-{n-k-1\choose k-1}+1\Big) \le sk{n-2 \choose k-2} \leq \frac{sk^3}{n^2}{n \choose k}.
        \end{equation}
We also have
\begin{equation}\label{eqhm2}{n\choose k}-{n-s\choose k} = {n\choose k}\Big(1-\prod_{i=0}^{k-1}\frac{n-s-i}{n-i}\Big)\ge {n\choose k}\Big(1-e^{-\sum_{i=0}^{k-1}\frac s{n-i}}\Big)\ge {n\choose k}\Big(1-e^{-\frac {sk}{n}}\Big).\end{equation}
The RHS of \eqref{eqhm1} is smaller than the RHS of \eqref{eqhm2} iff $1-\frac{sk^3}{n^2}>e^{-sk/n}.$
We have
\begin{equation} \label{lhs}
    1 - \frac{sk^3}{n^2} = \Big(1+\frac{sk^3}{n^2-sk^3}\Big)^{-1} \ge  e^{\frac{-sk^3}{n^2 - sk^3}}.
\end{equation}
Thus, the conclusion of the proposition holds, as long as
$$\frac{sk^3}{n^2 - sk^3}< \frac{sk}n$$
We solve the quadratic equation $n^2 - nk^2 - sk^3 > 0$ and get $n > \frac{k^2}{2}(1 + \sqrt{1 + \frac{4s}{k}})$. For $n > k^2 + k^{3/2}\sqrt{s}$ this inequality holds.
\end{proof}
We call families $\mathcal{A}, \mathcal{B} \subset 2^X$ \emph{cross-intersecting}, if $A \cap B \neq \emptyset$ for any $A \in \mathcal{A}, B \in \mathcal{B}$. We say that $C \subset X$ \emph{covers} family $\FF \subset 2^X$ if $C \cap A \neq \emptyset$ for any $A \in \FF$. We denote the size of the smallest cover of $\FF$ by $\tau(\FF)$.
We say that a family $\m{H}$ \emph{set-covers} a family $\m{F}$ if for each $F \in \m{F}$ there exists $H \in \m{H}$ such that $H \subset F$. For a graph $G$ denote $\mathcal I_k(G)$ the family of  independent $k$-element sets of vertices of $G$. In the proofs of Theorems~\ref{main} and~\ref{main2}, we need the following  result from \cite{KK}:
\begin{theorem}[Lemma~12 in \cite{KK}]\label{part}
Let $\m{F} \subset {[n] \choose k}$ be an intersecting family with $\tau(\m{F}) \geq 2$. Then we can split it into $\m{F}^{'} \sqcup \m{F}^{''}$, where $\m{F}^{'}$ can be set-covered by at most $k$ $2$-edges and $\m{F}^{''}$ cross-intersects some family of $t$-sets $\G^{\times}$ with $t \in \{ k-1, k\}$ and $\tau(\G^{\times}) \geq \sqrt{k}$. 
\end{theorem}

We can deduce the following structural property for the family $\ff''$ from the theorem above.

\begin{lemma}
\label{setcover}
    Let $\ff$ be a family of $k$-sets, which cross-intersects a family of $t$-sets $\G^{\times}$ with $t \in \{ k-1, k\}$ and $\tau({\G^{\times}}) \geq \sqrt{k}$. Then $\ff$ can be set-covered by a family $\m{H}$ of at most $k^{\sqrt{k}}$ $\sqrt{k}$-sets. 
\end{lemma}
\begin{proofw}
    
    We construct $\m H$ inductively. Put $\m H_0 = \emptyset$. For each $l \in [0, \sqrt{k} - 1]$, we construct a set-cover $\m H_{l+1}\subset{[n]\choose l+1}$ for $\mathcal F$ from the set-cover $\m H_{l}\subset {[n]\choose l}$ that we assume to have constructed. We also assume that  $|\m H_l| \leq k^{l}$.  %Assume that we constructed a family $\G_{l}$ which set-covers $\ff$: 
    Put $\ff[H] = \{F\in \ff: H\subset F\}$. Using this notation, the family $\m H_l$ being a set-cover means that $\ff \subset \bigcup_{H \in \m H_l} \ff[H]$. Thus, it should be clear that $\m H_0$ satisfies both the condition on size and on being a set-cover.  For each $H \in \m H_l$ consider a set $G_{H} \in \G^{\times}$ disjoint with $H$. Such set exists since $|H| = l < \sqrt{k}$ and $\tau(\G^{\times}) \geq \sqrt{k}$. Put $\m H_{l+1} = \{ H \cup x : H \in \m H_{l}, x \in G_{H} \}$. Then it is easy to see that $|\m H_{l+1}| \leq t |\m H_l| \le k^{l+1}$ and 
    $$
        \ff \subset \bigcup_{H \in \m H_{l+1}} \ff[H].
    $$
    Finally, we put $\m{H} = \m H_{\sqrt{k}}$.
\end{proofw}

Combining two previous statements, we get the following corollary.
\begin{corollary} \label{decompforchoice}
 \label{cont}
    Let $\ff$ be an intersecting family of $k$-sets.  If $\ff$ is non-trivial then it is a subfamily of  $\m{A} \cup \m{B}$, where the family $\m{A}$ is the family of all $k$-sets that are set-covered by a family $\m{E} \subset {[n] \choose 2}$ with $|\m{E}| = k,$ %consists of all $k$-sets containing one of $2$-sets from a family $$ 
    and the family $\m{B}$ is the family of $k$-sets that are set-covered by a family $\m{H} \subset {[n] \choose \sqrt{k}}$ with $|\m{H}| = k^{\sqrt{k}}.$   If $\ff$ is trivial it is a subfamily of some star $\m{S}_x$. 
\end{corollary}

\subsection{Proof of Theorem~\ref{main}}
In the proof of the theorem, we shall need the following result from \cite{KK}.
\begin{theorem}[Theorem~13 in \cite{KK}]\label{kk5}
Let $n, k$ be sufficiently large integers such that $n\ge \big(2+\frac 5{\ln k}\big)k^2$ and let $s < n$.  Let $G$ be a graph on $[n]$ with at most $sk$  edges. 
% and let $\mathcal{I}_k(G)$ be the family of all $k$-sets in $G$ that are independent in $G$.
Let $\mathcal{H} \subset \mathcal I_k(G)$
be a family of independent $k$-sets in $G$, which cross-intersects some family $\mathcal{G}$ of $t$-sets with $\tau(\mathcal{G}) \geq \sqrt{k}$ and $t \leq k$. Then $|\mathcal{I}_k(G)| > e^{k^{0.1}/2}|\mathcal H|$.
\end{theorem}

We also need the following result from  by Khadžiivanov and Nikiforov \cite{KhN}:
\begin{theorem}\label{KN}[Khadžiivanov--Nikiforov]
For a given graph $G$ let $\gamma$ be the density $\frac{|E(G)|}{|V(G)|^2}$ and $\m{N}_k(G)$ be the number of cliques on $k$ vertices. Then, if $\gamma \geq \frac{k-2}{2(k-1)}$, we have
$$
    \m{N}_k(G) \geq \frac{1+(2\gamma -1) (k-1)}{k}\cdot|V(G)|\cdot \m{N}_{k-1}(G) \text{ and } \m{N}_{k-1}(G) > 0.
$$
\end{theorem}

We prove Theorem~\ref{main} in assumption that $n < 2k^3$. The case $n \geq 2k^3$ is covered by  Proposition \href{hmbound}{\ref{hmbound}}.

\begin{proofw}
\underline{\emph{Proof of Theorem \href{main}{\ref{main}}}}.
We assume that $n < 2k^3$ and thus that $k$ itself is sufficiently large.
Let $\FF_1, \dots, \FF_s$ be a collection of $s$ non-trivial intersecting families. By Theorem \href{part}{\ref{part}}, we can split each $\FF_i$ into $\FF^{'}_i \sqcup \FF^{''}_i$, where $\FF^{'}_i$ can be set-covered by at most $k$ $2$-edges and $\FF^{''}_i$ cross-intersects some family of $t$-sets $\m{G}$ with $t \in \{ k-1, k\}$ and $\tau(\m{G}) \geq \sqrt{k}$. We get a split $\cup_{i \in [s]} \FF_i=\cup_{i \in [s]} (\FF^{'}_i \sqcup \FF^{''}_i)$. Next, we construct graph $G = ([n], E), |E|\le ks$, such that its set of edges $\mathcal{E}$ set-covers $\cup_{i \in [s]} \FF^{'}_i$. We denote by $\mathcal{I}_k(G)$  the family of independent $k$-sets of  $G$. %Consider an arbitrary set $A$ from $\m{I}_k(G)$. We may assume that $A$ is set-covered by an edge of $G$. Otherwise, it is already included in $\cup_{i \in [s]} \FF^{'}_i$. So 
We may of course assume that %$\cup_{i \in [s]} \FF^{''}_i$ consists of independent $k$-sets from $\m{I}_k(G)$, in other words 
$\cup_{i \in [s]} \FF^{''}_i \subset \m{I}_k(G) $.

By Theorem \href{kk5}{\ref{kk5}}, for any $\epsilon>0$, sufficiently large $k$ and $n>(2+\epsilon)k^2$ we have that  $|\mathcal I_k(G)| > k^4 |\FF^{''}_i|$ for any $i \in [s]$. Hence,
\begin{equation} \label{1}
\frac{s}{k^4}|\mathcal I_k(G)| > |\cup_{i \in [s]} \FF^{''}_i|.
\end{equation}
We need the following theorem, which we prove below.
\begin{theorem}\label{mainlemma}
Let $n, k$ be sufficiently large integers such that $2k^2< n < 2k^3$. Let $s < n$. % and let $n \geq \frac{1}{2}k + \frac 32 k \sqrt{s}$, $n>2k^2$.
Let $G$ be a graph on $[n]$ with at most $sk$ edges and let $\m{I}_k(G)$ be the family of all $k$-sets in $G$ that are independent in $G$.
 Then $|\mathcal I_k(G)| > \frac{k^4}{k^4 - s} {n-s \choose k}$.
\end{theorem}

Using Theorem \ref{mainlemma} and inequality \eqref{1}, it is straightforward to finish the proof of Theorem~\ref{main}.
$$
   \Big|\bigcup_{i \in [s]} \FF^{}_i\Big| = {n \choose k} - |\mathcal I_k(G)| +     \Big|\bigcup_{i \in [s]} \FF^{''}_i\Big| <  {n \choose k} -  \frac{k^4-s}{k^4-s}{n-s \choose k} = {n \choose k} - {n-s \choose k}.$$
\end{proofw}

Next, we prove Theorem \href{mainlemma}{\ref{mainlemma}}.
\begin{proofw}
\underline{\emph{Proof of Theorem \href{mainlemma}{\ref{mainlemma}}}}.
The quantity $|\m{I}_k(G)|$ is equal to  $\m{N}_k(\bar{G})$, where $\bar{G}$ is the complement of $G$. We iteratively apply %Khadžiivanov and Nikiforov's
Theorem \href{KN}{\ref{KN}}  in order to bound $\mathcal N_k(\bar{G})$:
$$
\m{N}_k(\bar{G}) \geq \frac{1+(2\gamma-1)(k-1)}{k}|V(\bar{G})|\m{N}_{k-1}(\bar{G}),
$$
$$
\dots
$$
$$
\m{N}_3(\bar{G}) \geq \frac{1+(2\gamma-1)(3-1)}{3}|V(\bar{G})|\m{N}_{2}(\bar{G}).
$$
In order to apply Theorem \href{KN}{\ref{KN}}, we need the inequality $\gamma \geq \frac{k-2}{2(k-1)}$ to hold. We have $\gamma = ({n\choose 2}-ks)/n^2$. Rearranging, we get a quadratic inequality $n^2-(k-1)n-2k(k-1)s\ge 0$. Since $s<n$, we have $n^2-(k-1)n-2k(k-1)s> n(n-2k(k-1))-(k-1)n=n(n-(2k+1)(k-1))>0$ for $n>2k^2$.  %It is valid for $n \geq \frac k 2+\frac 32 k\sqrt{s}$. 
In what follows, we use $\m{N}_{2}(\bar{G}) = |E(\bar{G})|= {{n \choose 2} - ks}$ and $2\gamma -1 = -\frac 1n-\frac{2ks}{n^2}$.
\begin{align*}
    \m{N}_k(\bar{G}) \geq& n^{k-2}\left( {n \choose 2} - ks \right)\prod_{i=2}^{k-1}\frac{1+(2\gamma -1)i}{i+1}\\
    % =& (n^k-n^{k-1}-2ksn^{k-2})\prod_{i=2}^{k-1}\left( i\left( 1-\frac{1}{n}-\frac{2ks}{n^2} \right) - i + 1 \right) \frac{1}{k!}\\
    =& \frac{n^k}{k!}\cdot \prod_{i=1}^{k-1}\left(1- i\left(\frac{1}{n}+\frac{2ks}{n^2} \right) \right).
\end{align*}

In what follows, we assume that ${n-s\choose k}>0$ (equivalent to $n\ge s+k$), otherwise the statement is trivial. In order to prove the theorem, it is thus sufficient for us to show that
$$
    \prod_{i=1}^{k-1}\left(1- i\left(\frac{1}{n}+\frac{2ks}{n^2} \right) \right) \cdot \left(1 - \frac{s}{k^4}\right) > \frac{k!{n-s \choose k}}{n^k} = \prod_{i=0}^{k-1}\left(1 - \frac{s+i}{n}\right).
$$
Since $n<2k^3\le k^4$, we have $\big(1-\frac s{k^4}\big)>\big(1-\frac sn\big)$, and thus the above inequality is implied by 
$$
    \prod_{i=1}^{k-1}\left(1- i\left(\frac{1}{n}+\frac{2ks}{n^2} \right) \right) > \prod_{i=1}^{k-1}\left(1 - \frac{s+i}{n}\right).
$$
We will simply show this term by term: for each $i\in[k-1]$, we have
$$ 1- i\left(\frac{1}{n}+\frac{2ks}{n^2} \right) > 1 - \frac{s+i}{n}.
$$
This is equivalent to $sn>2kis$, or $n>2ki$, which is valid since $n>2k^2\ge 2ki$.
\end{proofw}

\subsection{Proof of Theorem~\ref{main3}}
In order to prove the theorem, we give an explicit construction of a collection of non-trivial intersecting families whose union is bigger than that of the same number of trivial intersecting families. 

\begin{figure}
    \centering
\scalebox{.8}{
\begin{tikzpicture}
\coordinate (S) at (2, 2) ;
\coordinate (A) at (3.6, 0) ;
\coordinate (B) at (7.2,0) ;
\draw[thick, red] (0,0) rectangle (S);
\draw (1, 2) node[above] {$T_1$};

\fill[black] (0.4, 0.4)++(B) circle (0.1);
\fill[black] (1.6, 0.4)++(B) circle (0.1);
\fill[black] (1, 1.6)++(B) circle (0.1);

\fill[black] (0.4, 0.4)++(A) circle (0.1);
\fill[black] (1.6, 0.4)++(A) circle (0.1);
\fill[black] (1, 1.6)++(A) circle (0.1);

\fill[black] (0.4, 0.4) circle (0.1);
\fill[black] (1.6, 0.4) circle (0.1);
\fill[black] (1, 1.6) circle (0.1);

\node[ellipse, draw, color = green, thick, minimum width = 2cm, 
	minimum height = 0.8cm, rotate=65] at (4.3, 1) {};
\node[ellipse, draw, color = green, thick, minimum width = 2cm, 
	minimum height = 0.8cm, rotate=115] (e) at (4.9, 1) {};
\node[ellipse, draw, color = green, thick, minimum width = 1.8cm, 
	minimum height = 0.8cm] (e) at (4.58, 0.5) {};

\draw foreach \x in {1, 2, 3}
{ (2, 0)++(0.4*\x, 0) node {$\bullet$}};

\draw[thick, red] (A) rectangle ++(S);
\draw (A)++(1, 2) node[above] {$T_i$};

\draw foreach \x in {1, 2, 3}
{ (A)++(2 + 0.4*\x, 0) node {$\bullet$}};

\draw[thick, red] (B) rectangle ++(S);
\draw (B)++(1, 2) node[above] {$T_{k-1}$};

\coordinate (S`) at (2, -4);

\draw[thick, red] (0,-1) rectangle ++(S`);
\draw (1, -1) node[above] {$B_1$};

\draw foreach \x in {1, 2, 3}
{ (2, -5)++(0.4*\x, 0) node {$\bullet$}};

\draw[thick, red] (3.6, -1) rectangle ++(S`);
\draw (4.6, -1) node[above] {$B_i$};

\fill[black] (4.6, -1.3) circle (0.1);
\fill[blue] (4.6, -1.8) circle (0.1);
\draw[blue] (4.6, -1.8) node[right] {$y$};
\fill[black] (4.6, -2.3) circle (0.1);
\node at (4.6, -2.7) {\vdots};
\node at (4.6, -3.1) {\vdots};
\node at (4.6, -3.5) {\vdots};
\fill[black] (4.6, -4) circle (0.1);
\fill[black] (4.6, -4.5) circle (0.1);

\fill[black] (1, -1.3) circle (0.1);
\fill[black] (1, -1.8) circle (0.1);
\fill[black] (1, -2.3) circle (0.1);
\node at (1, -2.7) {\vdots};
\node at (1, -3.1) {\vdots};
\node at (1, -3.5) {\vdots};
\fill[black] (1, -4) circle (0.1);
\fill[black] (1, -4.5) circle (0.1);

\fill[black] (8.2, -1.3) circle (0.1);
\fill[black] (8.2, -1.8) circle (0.1);
\fill[black] (8.2, -2.3) circle (0.1);
\node at (8.2, -2.7) {\vdots};
\node at (8.2, -3.1) {\vdots};
\node at (8.2, -3.5) {\vdots};
\fill[black] (8.2, -4) circle (0.1);
\fill[black] (8.2, -4.5) circle (0.1);

\draw foreach \x in {1, 2, 3}
{ (5.6, -5)++(0.4*\x, 0) node {$\bullet$}};

\draw[thick, red] (7.2, -1) rectangle ++(S`);
\draw (8.2, -1) node[above] {$B_{k-1}$};

\draw[thick, red] (7.2, -6) rectangle ++(-5, -2);
\draw (4.6, -6) node[above] {$M$};

\draw [decorate,
	decoration = {brace, mirror}, thick] (10,0) --  (10,2);
\draw (10, 1) node[right] {$3$};
\draw [decorate,
	decoration = {brace}, thick] (10,-1) --  (10,-5);
\draw (10, -3) node[right] {$k-2-m$};
\draw [decorate,
	decoration = {brace}, thick] (10,-6) --  (10,-8);
\draw (10, -6) node[right] {$m$};

\draw [very thick, rounded corners, blue] (7.4, -5.7)--(7.4, -8.3)--(2,-8.3)--(2, -5.7)--(3.4, -5.7)--(3.4, -2)--(5, -2)--(5, -1.5)--(4.3, -1.5)--(3.4, -1.5)--(3.4, 2.6)--(5.8, 2.6)--(5.8, -5.7)--cycle;
\draw[ultra thick, blue] (5.9, -3) node[right] {$\boldsymbol{H_y^{i}}$};

\fill[black] (2.6, -7) circle (0.1);
\fill[black] (3.1, -7) circle (0.1);
\fill[black] (3.6, -7) circle (0.1);
\node at (4.1, -7) {\dots};
\node at (4.6, -7) {\dots};
\node at (5.1, -7) {\dots};
\node at (5.6, -7) {\dots};
\fill[black] (6.1, -7) circle (0.1);
\fill[black] (6.6, -7) circle (0.1);
 \end{tikzpicture}
 }
    \caption{Construction of intersecting families whose union is bigger than the union of stars. Red boxes correspond to the parts of the ground set. Their sizes are indicated on the right. Green ellipses are the 2-element sets that define $\mathcal T_i$. Blue vertex $y$ and a set $H_y^i$ define $\mathcal{HM}_y^i$.}
    \label{fig1}
\end{figure}
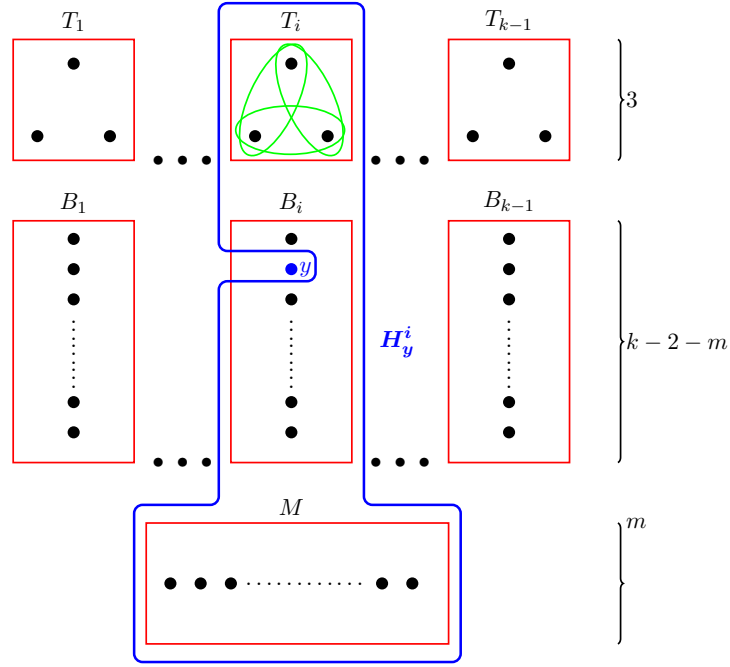

Let $n,k$ be positive integers, $m$ a non-negative integer, $m\le k-2$ and $n = (k-1)(k+1-m)+m$. Take the following pairwise disjoint sets: sets $T_1,\ldots, T_{k-1}$ of size $3$; sets $B_1,\ldots, B_{k-1}$ of size $k-2-m$; a set $M$ of size $m$. The ground set for the construction is the union of all these sets: $$G = M\sqcup \bigsqcup_{i=1}^{k-1}(T_i\sqcup B_i).$$
Note that $|G| = n$. The majority of the families in the construction would be of Hilton-Milner type. In order to define them, we need to define the `exceptional' sets. For  each $i\in [k-1]$ and $y\in B_i$ put $$H^i_y:= M\cup T_i\cup B_i\setminus \{y\}.$$
Consider the following collection of intersecting families (cf. Figure~\ref{fig1}). For each $i\in [ k-1]$ and $y\in B_i$, define
$$\mathcal{HM}^i_y = \{H_y^i\}\cup\Big\{A\in {G\choose k}: y\in A, A\cap H_y^i\ne \emptyset\Big\}.$$
Note that there are $\sum_{i=1}^{k-1}|B_i| = (k-1)(k-2-m)$ such families. Next, for each $i=1,\ldots, k-1$ define a family $\mathcal T_i$ of the form
$$\mathcal T_i:= \Big\{A\in {G\choose k}: |A\cap T_i| \ge 2\Big\}.$$

It is not difficult to see that the families $\mathcal{HM}_y^i$ and $\mathcal T_i$ are all non-trivial intersecting, moreover, their total number is $(k-1)(k-2-m)+k-1 = n-2k+2-m.$ 

Let us determine, how many sets $A\in {G\choose k}$ are not covered by these  families. If $|A\cap T_i|\ge 2$ then $A$ is covered by $\mathcal T_i$. $y\in A\cap B_i$ and $|A\cap (T_i\sqcup B_i)|\ge 2$ then $A \in \mathcal{HM}_y^i$. Thus, all $A$ with  $|A\cap (T_i\sqcup B_i)|\ge 2$ are covered by our families.

If  $A$ is such that  $|A\cap (T_i\sqcup B_i)|\le 1$  for all $i\in [k-1]$, then, by the pigeon-hole principle, $|A\cap M|\ge 1$. If at the same time $y\in A\cap B_i$ for some $i$ then $A \in \mathcal{HM}_y^i$.

Concluding, $A$ is not covered by these families if (and only if) $A\subset M\sqcup \bigsqcup_{i=1}^{k-1}T_i$ and $|A\cap T_i|\le 1$ for each $i\in[k-1]$. The number of such sets is $$\sum_{j=1}^m{m\choose j}{k-1\choose k-j}3^{k-j}.$$ For $m=1$ this sum is equal to $3^{k-1}$. On the other hand, the number of sets not covered by a union of $n-2k+2-m$ stars is $${2k-2+m\choose k}.$$ For $m=1$, it is equal to $(4+o(1))^k$. It is not difficult, although tedious, to check that the first expression is smaller than the second for large enough $k$ and $m\le k/4$. We decided to omit these calculations.

\subsection{Proof of Theorem~\ref{main2}}
 We employ the probabilistic method. Put $u := Cn\log{n}$ with some $C$ that we specify later on. Take a set of $u$ colors and correspond to each vertex of $KG_{n, k}$ a random subset of colors of size $\frac{u}{2}$. Any independent set of vertices of $KG_{n, k}$ is an intersecting family in ${[n] \choose k}$. Using Corollary \ref{decompforchoice}, we get that any such independent set is contained in one of the families from $\m{C}$, where $\m{C}$ consists of families of the form $\m{A} \cup \m{B}$ or $\m{S}_x$.
 Note that 
 $$
    |\m{C}| \leq {{n \choose 2} \choose k}{{n \choose \sqrt{k}} \choose k^{\sqrt{k}}} + n. 
$$
We say that a partition $X = X_1 \sqcup \dots \sqcup X_m$ of a set $X$ {\it lies in a cover} of this set $X = X'_1 \cup \dots \cup X'_m$, if 
$X_i \subset X'_i$ for each $i\in[m]$. Hence, any partition of ${[n] \choose k}$ lies in one of  
$$
    |\m{C}|^u \leq \left({{n \choose 2} \choose k}{{n \choose \sqrt{k}} \choose k^{\sqrt{k}}} + n\right)^u \leq n^{\left(2k + k^{\sqrt{k} + \frac{1}{2}}\right)u} %= n^{k^{\sqrt{k}(1 + o(1))}}
$$
possible covers, formed by a $u$-tuple of families from $\m{C}$.

For a fixed cover $\m{K} = \m{K}_1 \cup \dots \cup \m{K}_u$ from $\m{C}$ we are going to bound from above the probability of the event $\m{A_K}$ that $KG_{n,k}$ can be colored in one of the colorings of ${[n] \choose k}$ that lies in $\m{K}$ and that respects the lists of colors assigned to the $k$-sets. 

First, we show that families of type $\m{B}$ cover only a small part of vertices of $KG_{n, k}$. Consider an arbitrary family of type $\m{B}$. Clearly, it covers less than $k^{\sqrt{k}}{n-\sqrt{k} \choose k - \sqrt{k}}$ $k$-sets. We have at most $u$ families of such type, and thus together they can cover at most 
$$
    uk^{\sqrt{k}}{n-\sqrt{k} \choose k - \sqrt{k}} \leq Cn\log{n} \: k^{\sqrt{k}} \left(\frac{k}{n}\right)^{\sqrt{k}}{n \choose k}
$$
sets. We have 
$$
    \frac{Cn\log{n} \: k^{2\sqrt{k}}}{n^{\sqrt{k}}} <  \frac 16
$$
for large enough $n$ that also satisfies $n > (1 + \varepsilon)k^2$. Therefore, the number of $k$-sets covered by families of type $\m{B}$ is less than $\frac 16{n \choose k}$.

Second, we analyze, what could be covered by the families  $\m{A}$ and $\m{S}_x$. Consider a random $k$-set $S$ from $[n]$. Let $\eta_1$ denote the number of $x\in S$ such that $x$ is a center of some 
$\mathcal S_x$. Let $\eta_2$ denote the number of pairs $\{x,y\}\subset S$ such that $\{x,y\}$ is contained in $\m{E}(\m{A})$ for some $\m{A}.$  For convenience, we also call such pairs \emph{centers}. Then the random variable $\eta_1+\eta_2$ equals to the number of \emph{centers} of type $\mathcal S_x$ or $\m{A}$ in a random $k$-set.

Suppose that there are $u_1$ families of type $\mathcal S_x$ and $u_2$ families of type $\mathcal A$ in the cover $\m K$. Using that each family of type $\m A$ consists of at most $k$ pairs and that $n\ge k^2$, we can upper bound the expectation of $\eta_1+\eta_2$:
$$\E[\eta_1+\eta_2] \le u_1 \frac kn+u_2 \frac {k{k\choose 2}}{{n\choose 2}}\le  u_1 \frac kn+u_2 \frac {k^3}{n^2}\le u_1 k/n+u_2 k/n\le u k/n.$$
Using Markov’s inequality, we get that
$$
\Pr\left[\text{ a uniformly random set from } {[n] \choose k} \text{ contains more than } \frac{2uk}{n} \text{ centers }\right] \leq \frac{1}{2}.
$$
Thus, at least half of all ${[n] \choose k}$ sets contain at most $\frac{2uk}{n}$ \emph{centers}. It means that each of such sets can be colored in at most $\frac{2uk}{n}$ colors of type $\mathcal S_x$ or $\m{A}$. (At this point, we ignore the color lists assigned to the vertices.) From the previous paragraph, we also know that most  of them cannot be colored in colors of type $\m{B}$. Let $\m{X}$  denote the family of sets which can only be colored in at most $\frac{2uk}{n}$ colors of types $\mathcal S_x$ or $\m{A}$. Then 
$$
   |\m{X}| \geq \frac{1}{3}  {n \choose k}.
$$ 

Given $X \in \m{X}$ and the color list assigned to $X$, let $B_{X, \m{K}}$  stand for the event that $X$ cannot be colored with a coloring that lies in $\m{K}$. $X$ can be colored in color $j$, where $\mathcal K_j$ of type $\mathcal S_x$ or $\mathcal A$ only if $j$ belongs to the randomly chosen subset of $\frac u 2$ colors that was assigned to $\m{X}$. There are $2u k/n =  2Ck\log_2n=:z$ colors that contain $X$. Then the probability $p$ of the event $B_{X, \m{K}}$ is at least the probability that the colors for $X$ were chosen from the
complement of the set of possible colors:
 $$p\ge \frac{{u- z\choose u/2}}{{u\choose u/2}} = 2^{-z}(1-o(1)) \geq \frac 12 n^{-2Ck }.$$

Recall that  $\m{A_K}$ is the event that $KG_{n,k}$ can be colored in one of the colorings of ${[n] \choose k}$ that lies in $\m{K}$ and that respects the lists of colors assigned to the $k$-sets. We have $\m{A}_{\m{K}} \subset \bigcap_{X \in \m{X}}\bar{B}_{X, \m{K}} $, and thus
 $$\Pr [\m{A}_{\m{K}}] \leq (1-p)^{|\m{X}|} \leq e^{-p|\m{X}|}.$$

We have 
 $$
    p|\m{X}| \geq \frac 16 n^{-2Ck }  {n \choose k} \geq \frac 16 n^{-2Ck} \left(\frac{n}{k}\right)^k \geq n^{k/2 -2Ck},
 $$
for large enough $n > (1+\varepsilon)k^2$, therefore
$$
    \Pr [\m{A}_{\m{K}}] \leq e^{-n^{k/2 -2Ck}}.
$$

Let $U$ denote the event that there exists a proper coloring of $KG_{n,k}$ that respects the lists assigned to the vertices. 

To conclude the proof, it is enough to verify the last inequality below.  
$$
    \Pr[U] \leq \sum_{\m{K} \in \m{C}^u} \Pr [\m{A}_{\m{K}}] \leq n^{\left(2k + k^{\sqrt{k} + \frac{1}{2}}\right)u} \cdot e^{-n^{k/2 -2Ck}} < 1.
$$
Comparing the exponents, it is equivalent to
$$2k + k^{\sqrt{k}+1/2} < \frac{n^{k/2 - 1 - 2Ck}}{C\log^2{n}}.
$$
The RHS is at least $2n^{k(1/2  - 2C) - 1 - \delta}$ for an arbitrary $\delta>0$ and $n\ge n_0(\delta)$, and the LHS is at most $2 k^{\sqrt{k}+1/2}$. Thus, it is sufficient to check the inequality
\begin{equation} \label{eqfin}
   k^{\sqrt{k}+1/2}  < n^{k(1/2  - 2C) - 1 - \delta}.
\end{equation}
It is easy to check that it holds if $1/2-2C\ge \delta+\frac 1k$ and $n$ sufficiently large. The inequality we get on $C$ is $C\le \frac 14 - \frac 1{2k}-\epsilon.$ For $k=3$, for example, this gives $\frac 1{12}-\epsilon.$ The theorem is proved.

%\printbibliography

\end{document}